\documentclass{article}

\usepackage{amsmath,amsthm}
\usepackage{caption2}
\usepackage[dvips]{color}
\usepackage{amsfonts}
\usepackage{rotating}
\usepackage{euscript}
\usepackage{pst-node}
\usepackage{epsfig}
\usepackage{comment}

\def\be{\begin{equation}}
\def\ee{\end{equation}}

\def\d{{d\,}}

\def\C{{\mathbb C}}

\def\P{{\mathbb P}}
\def\Z{{\mathbb Z}}

\def\Q{{\mathbb Q}}

\def\phi{{\varphi}}

\def\deg{{\rm deg\,}}

\def\dim{{\rm dim\,}}

\def\bp{\begin{proposition}}
\def\ep{\end{proposition}}

\def\bt{\begin{theorem}}
\def\et{\end{theorem}}
\def\be{\begin{equation}}
\def\bee{\begin{equation*}}
\def\la{\label}

\def\ee{\end{equation}}
\def\eee{\end{equation*}}
\def\bl{\begin{lemma}}
\def\el{\end{lemma}}
\def\bc{\begin{corollary}}
\def\ec{\end{corollary}}

\def\bd{\begin{definition}}
\def\ed{\end{definition}}
\def\wt{\widetilde}
\newcommand{\vv}{\vec{v}}
\newcommand{\vw}{{\vec{w}}}

\input epsf.sty

\newtheorem{theorem}{Theorem}[section]
\newtheorem{lemma}[theorem]{Lemma}
\newtheorem{definition}[theorem]{Definition}
\newtheorem{corollary}[theorem]{Corollary}
\newtheorem{proposition}[theorem]{Proposition}

\newcommand{\ve}{\vec{e}}

\title{Laurent polynomial moment problem: \\ a case study}

\author{
F.\,Pakovich\thanks{Faculty of Natural Sciences, Ben-Gurion University
of the Negev, P.O.B. 653, Beer-Sheva, Israel; e-mail:
{\tt pakovich@cs.bgu.ac.il}}, 
C.\,Pech\thanks{Department of Algebra, Johannes Kepler University, 
Altenberger Strasse 69, 4040, Linz, Austria,
e-mail: {\tt cpech@freenet.de}}, 
A.\,Zvonkin\thanks{LaBRI, Universit\'e Bordeaux I, 351 cours de la 
Lib\'eration, F-33405 Talence Cedex France; e-mail: 
{\tt zvonkin@labri.fr}}
}

\date{}

\begin{document}

\maketitle

\begin{abstract}
In recent years, the so-called {\em polynomial moment problem}, motivated
by the classical Poincar\'e center-focus problem, was thoroughly studied,
and the answers to the main questions have been found. The study of a 
similar problem for rational functions is still at its very beginning. 
In this paper, we make certain progress in this direction; namely, we 
construct an example of a Laurent polynomial for which the solutions of 
the corresponding moment problem behave in a significantly more 
complicated way than it would be possible for a polynomial.
\end{abstract}

\section{Introduction} 

The main result of this paper is a construction of a particular 
Laurent polynomial with certain unusual properties. This Laurent 
polynomial is a counter\-example to an idea that, 
so far as the moment problem is concerned, rational functions would 
behave in the same way as polynomials.
The main interest of the paper, besides the result itself, lies 
in a peculiar combination of methods which involve certainly the complex
functions theory but also group representations, Galois theory, 
and the theory of Belyi functions and ``dessins d'enfants'', 
while the motivation for the study comes from differential equations.

In addition to theoretical considerations our project involves 
computer calculations. It would be difficult to present here all 
the details. However, we tried to supply an interested reader with
sufficient number of indications in order for him or her 
to be able to reproduce our results. A less interested reader
may omit certain parts of the text and just take our word for it.

\pagebreak[4]

About a decade ago, M.\,Briskin, J.-P.\,Fran\c{c}oise, and Y.\,Yomdin 
in a series of papers \cite{bfy1}--\cite{bfy4} posed the following
 
\paragraph{Polynomial moment problem.} 
{\it For a given complex polynomial\/ $P$ and distinct complex numbers\/
$a,b$, describe all polynomials\/ $Q$ such that}
\be\la{111}
\int^b_a P^i\,dQ \,=\, 0
\ee
{\it for all integer\/ $i\geq 0$.}

\bigskip

The polynomial moment problem is closely related to the center problem 
for the Abel differential equation in the complex domain, which in its 
turn may be considered as a simplified version of the classical 
Poincar\'e center-focus problem for polynomial vector fields.
The center problem for the Abel equation and the polynomial moment 
problem have been studied in many recent papers (see, e.\,g., 
\cite{bby}--\cite{by}, \cite{c}, \cite{mp}--\cite{y1}).

There is a natural sufficient condition for a polynomial $Q$ to 
satisfy~\eqref{111}. Namely, suppose that there exist polynomials 
$\wt P$, $\wt Q$, and $W$ such that
\be\la{222}
P=\wt P\circ W, \qquad Q=\wt Q\circ W, \qquad \mbox{and} \qquad  W(a)=W(b),
\ee
where the symbol $\circ$ denotes a superposition of functions: 
$f_1\circ f_2=f_1(f_2)$. Then, after a change of variables $z\to W(z)$ 
the integrals in~\eqref{111} are transformed to the integrals 
\be\la{pop}
\int^{W(b)}_{W(a)} \wt P^i\, d\wt Q
\ee 
and therefore vanish since the polynomials $\wt P^i$ and $\wt Q$ are 
analytic functions in~$\C$ and the integration path in~\eqref{pop} 
is closed. A solution of~\eqref{111} for which~\eqref{222} holds 
is called {\it reducible}. For ``generic'' collections $P$, $a$, $b$ any 
solution of~\eqref{111} turns out to be reducible. For instance, 
this is true if $a$ and $b$ are not critical points of~$P$, see~\cite{c}, 
or if $P$ is indecomposable, that is, if it cannot be represented as
a superposition of two polynomials of degree greater than one,
see~\cite{pa2} (in this case (\ref{222}) reduces to the equalities
$P=W$, $Q=\wt Q \circ P$, and $P(a)=P(b)$).
Nevertheless, as it was shown in~\cite{pa1}, if $P(z)$ has several 
composition factors $W$ such that $W(a)=W(b)$ then the sum of the 
corresponding reducible solutions may be an irreducible one.

It was conjectured in \cite{pa4} that actually {\em any}\/ solution 
of~\eqref{111} can be represented as a sum of reducible ones. Recently
this conjecture was proved in~\cite{mp}. The proof relies on two key 
components. The first one is a result of~\cite{pa3} which states that 
$Q$ satisfies (\ref{111}) if and only if the superpositions of 
$Q$ with branches $P_i^{-1}(z)$, $1 \le i \le n$, of the algebraic 
function $P^{-1}(z)$ satisfy a certain system of linear equations 
\be \la{su}
\sum_{i=1}^nf_{s,i}Q(P^{-1}_{i}(z))=0, \ \ \ \ \ \ f_{s,i} \in \Z,  
\ \ \ \ \ \ 1\leq s\leq k, 
\ee  
associated to the triple $P$, $a$, $b$ in effective way.

The second key componenet is related to the vector subspace 
$V_{P,a,b}\subset \mathbb Q^n$ spanned by the vectors 
$$
(f_{s,\sigma(1)}, f_{s,\sigma(2)},\, ...\,, f_{s,\sigma(n)}), \qquad 
1\leq s \leq k, \qquad \sigma\in G_P\,,
$$ 
where $G_P$ is the monodromy group of $P$ and $f_s,$ $1\leq s \leq k,$ 
are vectors from~\eqref{su}.
By construction, the subspace $V_{P,a,b}$ is invariant under the action 
of $G_P$, so the idea is to obtain
a full description of such subspaces.
In short, it was proved in~\cite{mp} that if a transitive permutation 
group~$G\le {\rm S}_n$ contains a cycle of length~$n$ then the 
decomposition of $\mathbb Q^n$ in irreducible components of the 
action of $G$ depends only on the imprimitivity systems of $G$. 
Obviously, the monodromy group of a polynomial of degree $n$ always 
contains a cycle of length $n$ which corresponds to the loop around 
infinity. Furthermore, imprimitivity systems of $G_P$ correspond to 
functional decompositions of $P$. Therefore, the structure of invariant 
subspaces of the permutation representation of $G_P$ over $\Q$ depends 
only on the structure of functional decompositions of $P$, and a careful 
analysis of system~\eqref{su} and of the associated space $V_{P,a,b}$ 
eventually permits to prove that any solution of~\eqref{111} is a sum 
of reducible solutions.
Notice that using the decomposition theory of polynomials one can also 
show that actually any solution of~\eqref{111} may be represented as 
a sum of at most {\it two}\/ reducible solutions, and describe these 
solutions in a very explicit form (see \cite{rii}).

For example, in the simplest case of the problem corresponding to an 
{\it indecomposable} polynomial $P$ the above strategy works as follows. 
The only invariant subspaces of the permutation representation of $G_P$ 
on $\Q^n$ in this case are the subspace $V_1$ spanned by the vector 
$(1,1,\dots, 1)$, and its complement $V_1^{\bot}$. Since system~\eqref{su} 
contains an equation whose coefficients are not all equal this implies that
$V_{P,a,b}=V_1^{\bot}$ and therefore \eqref{su} yields that 
$$
Q(P^{-1}_{1}(z))= Q(P^{-1}_{2}(z))=\dots =Q(P^{-1}_{n}(z))
$$ 
identically over $z$. On the other hand, such an equality is possible 
only if $Q=\wt Q\circ P$ for 
some $\wt Q\in \C[z]$. Finally, $P(a)=P(b)$ since otherwise after the change 
of variables $z\rightarrow P(z)$ we would obtain that $\wt Q$ is orthogonal 
to all powers of $z$ on $[P(a),P(b)]$ in contradiction to the Weierstrass 
theorem.

In the paper \cite{pak2} the following generalization of the 
polynomial moment problem was 
investigated: for a given rational function $F$ and a curve 
$\gamma\subset \C\P^1$, describe rational functions $H$ such that 
\be\la{mom}
\int_{\gamma} F^i\,dH \,=\, 0
\ee 
for all $i\geq 0$. 
In particular, in \cite{pak2} another version of system \eqref{su} was 
constructed: its solutions, instead of the equality \eqref{mom},
guarantee only the rationality of the generating function 
$f(t)=\sum_{i=0}^{\infty}m_it^i$ for the moments 
\be\la{mom1}
m_i=\int_{\gamma} F^i\,dH \,.
\ee
On the other hand, it was shown that if the additional conditions 
$H^{-1}\{\infty\}\subseteq F^{-1}\{\infty\}$ and $F(\infty)=\infty$ are 
satisfied, then the rationality of $f(t)$ actually implies that 
$f(t)\equiv 0$.

The following modification of~\eqref{222} is a natural sufficient 
condition implying~\eqref{mom}: there exist rational functions 
$\wt F$, $\wt H$, and $W$ such that
\be\la{ge}
F \,=\, \wt F \circ W, \qquad H \,=\, \wt H \circ W\,,
\ee
the curve $W(\gamma)$ is closed, and all the poles of the functions 
$\wt F$, $\wt H$ lie ``outside''~$W(\gamma)$
(the term ``outside'' is written in quotation marks since it is
defined also for self-intersecting curves).
We will call such a solution of~\eqref{mom} {\it geometrically}\/ 
reducible. Notice that if $\gamma$ is closed then geometrically 
reducible solutions always exist. Indeed, one may take
$$
W \,=\, F, \qquad H \,=\, \wt H \circ F
$$ 
where $\wt H$ is any rational function with all its poles outside
the curve $F(\gamma)$. 
It is also shown in \cite{pak2} that, similarly to the case of a
polynomial~$P$, for a generic rational function~$F$ (for example, 
for an $F$ whose monodromy group is the full symmetric group) all 
solutions of~\eqref{mom} turn out to be geometrically reducible.
However, for a non-generic $F$ the situation becomes much more 
complicated in comparison with the polynomial moment problem, and some 
reasonable description of solutions of~\eqref{mom} seems (at least
for the moment) to be unachievable. 

\bigskip

In this paper we will consider a particular case of problem~\eqref{mom} 
which is especially interesting in view of its connection with the 
classical version of the Poincar\'e center-focus problem. Namely, 
we will consider the following 

\paragraph{Laurent polynomial moment problem.}
{\it For a given Laurent polynomial\/~$L$ which is not a polynomial in\/ 
$z$ or in\/ $1/z$, describe all Laurent polynomials\/~$Q$ such that}
\be\la{1l}
\int_{S^1} L^i\,dQ \,=\, 0
\ee
{\it for all integer\/ $i\geq 0$}.

\bigskip
 
In contrast to the polynomial moment problem, not any solution of the 
Laurent polynomial moment problem is a sum of geometrically reducible 
solutions. For example, as it was observed in~\cite{pak2}, if 
$L(z)=\wt L(z^d)$ for some $d>1$, then the residue calculation shows 
that condition \eqref{1l} is satisfied for any Laurent polynomial $Q$ 
containing no terms of degrees which are multiples of $d$.
We will call such a solution of the Laurent polynomial moment problem 
{\it algebraically}\/ reducible. Notice that, in distinction to
geometrically reducible solutions which always exist, algebraically 
reducible solutions exist only if $L$ is decomposable and has $z^d$ as 
its right composition factor.
One might think that any solution of the Laurent polynomial moment 
problem is a sum of geometrically and/or algebraically reducible 
solutions but, as we will see below, this is not the case 
either, although it seems that for a ``majority'' of Laurent 
polynomials $L$ this is the case.

It is natural to start the investigation of the Laurent polynomial moment 
problem by the study of the particular case where $L$ is indecomposable. 
At least, in this case there exist no algebraically reducible solutions. 
On the other hand, any geometrically reducible solution of~\eqref{1l} 
must have the form $Q = \wt Q \circ L$, where~$\wt Q$ is a rational 
function whose poles lie outside the curve $L(S^1)$. However, since~$Q$ 
is a Laurent polynomial, it is easy to see that in this case $\wt Q$ 
is necessarily a polynomial. Furthermore, a sum of geometrically reducible 
solutions has the form 
$$
\sum_{i} \wt Q_i\circ L \,=\, \left(\sum_{i} \wt Q_i\right)\circ L
$$
and hence is itself a geometrically reducible solution. Thus, ``expectable'' 
(and therefore not very interesting) solutions of the Laurent polynomial 
moment problem for indecomposable $L$ are of the form $Q=\wt Q(L),$ 
where $\wt Q$ is a polynomial. Any other solutions, when they exist,
are of great interest since they show that the situation is more
complicated than one might hope.

Let $L$ be an indecomposable Laurent polynomial of degree $n$, and let 
$G_L$ be its monodromy group. We will always assume that $L$ is {\it proper},
that is, it has poles both at zero and at infinity.
In this case the group $G_L$ contains a permutation with two cycles: this
permutation corresponds to the loop around infinity.
Furthermore, it follows from Theorem 4.5 of~\cite{pak2} that if the only 
invariant subspaces of the permutation representation of $G_L$ on 
$\Q^n$ are $V_1$ and $V_1^{\bot}$, then any solution of the Laurent 
polynomial moment problem for $L$ is geometrically reducible.
Therefore, if we want to find an example of a Laurent polynomial $L$
for which there exist solutions which are not geometrically reducible, 
we may use the following strategy:
\begin{itemize}
\item   First, find a permutation group $G$ of degree $n$ such that $G$ 
        would contain a permutation with two cycles, and the permutation 
        representation of~$G$ on~$\Q^n$ would have more than two invariant 
        subspaces. 
\item   Then, realize $G$ as the monodromy group of a Laurent 
        polynomial~$L$. 
\item   And, finally, prove somehow the existence of non-reducible
        solutions.
\end{itemize}

This program was started in \cite{pak2}. Namely, basing 
on Riemann's existence theorem it was shown that there exists a Laurent 
polynomial~$L$ of degree 10 such that its monodromy group is permutation 
isomorphic to the action of\/ ${\rm S}_5$ on two-element subsets of the 
set of\/ 5 points. The corresponding permutation action of\/ ${\rm S}_5$ 
on\/ $\Q^{10}$ has more than two invariant subspaces. Furthermore, 
proceeding from a general algebraic result of Girstmair~\cite{g} about 
linear relations between roots of algebraic equations  it was shown 
that there exists a rational function~$Q$ which is not a rational function 
in $L$ such that the generating function for the sequence of the moments  
$$
m_i=\int_{S^1} L^i\,dQ, \ \ \ i\geq 0
$$
is rational (see Sec.~8.3 of \cite{pak2}). However, the methods 
of~\cite{pak2} do not permit to find $L$ or $Q$ explicitly and tell 
us nothing about the structure of solutions of~\eqref{1l}.

In this paper we provide a detailed analysis of the above example with 
the emphasis on the two following questions of a general nature:
\begin{itemize}
\item   First, how to construct a Laurent polynomial $L$ starting from its 
        monodromy group\/~$G_L$?
\item   Second, how to describe solutions of~\eqref{1l} which are not 
        geometrically reducible ?
\end{itemize}
We answer all these questions for the particular Laurent polynomial~$L$
given below. Actually, we believe that our methods can be used in a more 
general situation too and can serve as a ``case study'' for further research
concerning the Laurent polynomial moment problem.

\bigskip

The main result of this paper is an actual calculation of an indecomposable 
Laurent polynomial $L$ such that the corresponding moment problem has 
non-reducible solutions, and a complete description
of these solutions. Namely, we show that for
\be\label{expr-for-L}
  L \,=\, \frac {K \left( z-1 \right) ^{6} \left( z-a \right) ^{3} \left( z-b
  \right) }{{z}^{5}}\,,
\ee
where
$$
 K \,=\, {\frac {11}{216}}+{\frac {5}{216}}\,\sqrt {5}, \qquad
 a \,=\, -\frac{3}{2}+\frac{1}{2}\,\sqrt {5}, \qquad
 b \,=\, \frac{7}{2}-\frac{3}{2}\,\sqrt {5}\,,
$$
there exist Laurent polynomials $Q_0=1$, $Q_1, Q_2, Q_3, Q_4$
(we compute them explicitly in Sec.~\ref{sec:proof}) such that the 
following statement holds:

\bt\la{main} 
A Laurent polynomial $Q$ is orthogonal to all powers of $L$ on $S^1$ 
if and only if $Q$ can be represented in the form
$$
Q \,=\, \sum_{j=0}^4 (R_j \circ L) \cdot Q_j
$$
for some polynomials $R_0, R_1, R_2, R_3, R_4$.
\et

In other words, solutions of the moment problem for $L$ form a 
5-dimensional module over the ring of polynomials in $L$ (while in a 
generic case such a module is one-dimensional and is therefore composed 
of polynomials in $L$ and of nothing else). The choice of the basis 
$Q_j$ is not unique, but once a basis is chosen the above representation 
of $Q$ becomes unique. 

\bigskip

The paper is organized as follows. In Sec.~2 we give a detailed 
description of the permutation action of ${\rm S}_5$ on $\Q^{10}$.  
In Sec.~3 we compute explicitly a Laurent polynomial $L$ whose monodromy 
group is permutation equivalent to this action. Finally, in Sec.~4 we 
determine the above mentioned Laurent polynomials $Q_j$ and prove 
Theorem~\ref{main}.

\bigskip

\paragraph{Acknowledgments.}
The first author is grateful to C.\,Christopher, J.\,L.\,Bra\-vo, and 
M.\,Muzychuk for valuable discussions, and to the Max-Planck-Institut 
f\"ur Mathematik for hospitality. The second author is indebted to 
the Center for Advanced Studies in Mathematics of the Ben-Gurion 
University for the financial support during his Post Docorate
at Ben-Gurion University. The first and the third authors wish to
thank the International Centre for Mathematical Sciences 
(Edinburgh) for the possibility to discuss a preliminary version 
of this paper.

\section{Permutation representation of ${\rm S}_5$ on $\Q^{10}$
with more than two invariant subspaces}\la{dim}

\begin{figure}[htbp]
\centerline{\epsffile{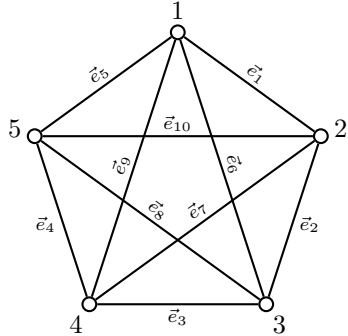}}
\caption{\small The correspondence between the edges of the complete 
graph $K_5$ and the basis in $\Q^{10}$.}
\label{corr}
\end{figure}

Consider the {\em complete graph}\/ $K_5=(V,E)$ having the vertex set
$V=\{1,2,3,4,5\}$ and the edge set $E$ consisting of all the subsets 
of $V$ of size~2. The symmetric
group~${\rm S}_5$ acts on $V$ and therefore also on $E$, and we thus 
obtain a transitive action of ${\rm S}_5$ of degree~10. Moreover,
the homomorphism ${\rm S}_5 \to {\rm S}_{10}$ is obviously injective.
Let us identify the canonical basis
\begin{eqnarray*}
 \ve_1 & = & (1,0,0,0,0,0,0,0,0,0)\,, \\
 \ve_2 & = & (0,1,0,0,0,0,0,0,0,0)\,, \\
 \vdots & \qquad & \vdots \\
 \ve_{10} & = & (0,0,0,0,0,0,0,0,0,1)
\end{eqnarray*}
of the space $\Q^{10}$ with the set $E$. This identification may in 
principle be arbitrary but we have chosen the one which is more 
``readable'', see~Fig.~\ref{corr}: the first five vectors are 
associated, in a cyclic way, to the sides of the pentagon, while 
the last five vectors are associated in the similar way to the sides 
of the inside pentagram.

Associating to each element of ${\rm S}_5$ a $10 \times 10$ permutation
matrix corresponding to the action of this element on $E$ we obtain a
{\em permutation representation}\/ of~${\rm S}_5$ on $\Q^{10}$.
Any permutation representation of any finite group always has at least 
two invariant subspaces: the subspace $U_1$ of dimension~1 spanned by 
the vector $\vec{\mathbf 1}=(1,1,\ldots,1)$, and its orthogonal 
complement $U_{n-1}=U_1^{\perp}$ of dimension~$n-1$ containing the vectors 
$(x_1,x_2,\ldots,x_n)$ having $\sum_{i=1}^nx_i=0$. While the space $U_1$ 
is obviously irreducible, the space $U_{n-1}$ may be, or may not be 
irreducible. We will show that in our case it is reducible.

\begin{figure}[htbp]
\centerline{\epsffile{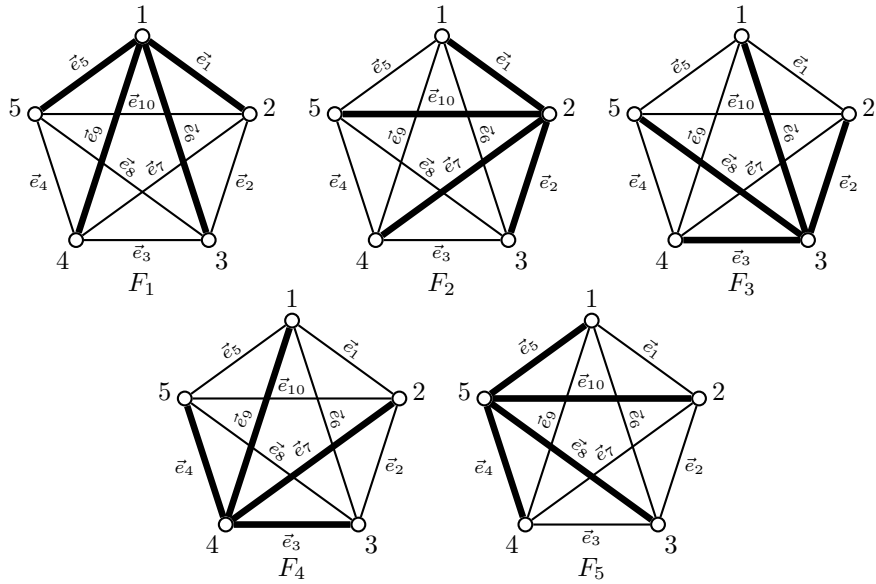}}
\caption{\small The fans $F_i$, that is, the sets of edges incident to 
the vertex $i=1,\ldots,5$.}
\label{corr2}
\end{figure}

One of the ways to construct invariant subspaces in our example is
to consider subsets of edges which are sent to one another by the
action of ${\rm S}_5$ on the vertices. Let us take the {\em fans}\/ 
$F_i \subset E$, $i=1,\ldots,5$, where $F_i$ is the set of edges of $K_5$
incident to the vertex~$i$, see Fig.~\ref{corr2}. 

Obviously, any permutation of the vertices sends fans to fans. Therefore, 
the vectors $\vv_i=\sum_{u\in F_i}e_u$, or, more concretely,
\begin{eqnarray*}
\vv_1 & = & (1,0,0,0,1,1,0,0,1,0)\,, \\
\vv_2 & = & (1,1,0,0,0,0,1,0,0,1)\,, \\
\vv_3 & = & (0,1,1,0,0,1,0,1,0,0)\,, \\
\vv_4 & = & (0,0,1,1,0,0,1,0,1,0)\,, \\
\vv_5 & = & (0,0,0,1,1,0,0,1,0,1)
\end{eqnarray*}
(the first five and the last five components of these vectors move 
cyclically) span an invariant subspace $F\subset\Q^{10}$. It is easy to 
verify that $F$ is 5-dimensional. Since every edge is contained in exactly 
two fans we have $\sum_{i=1}^5\vv_i=(2,2,\ldots,2)$ and therefore $F$ 
contains $U_1$ as its subspace. The orthogonal complement of $U_1$ in 
$F$ is a 4-dimensional invariant subspace $U_4\subset\Q^{10}$. 
The vectors
\begin{eqnarray*}
 \vv_2-\vv_1 & = & (0,1,0,0,-1,-1,1,0,-1,1)\,,\\
 \vv_3-\vv_1 & = & (-1,1,1,0,-1,0,0,1,-1,0)\,,\\
 \vv_4-\vv_1 & = & (-1,0,1,1,-1,-1,1,0,0,0)\,,\\
 \vv_5-\vv_1 & = & (-1,0,0,1,0,-1,0,1,-1,1)\,,
\end{eqnarray*}
each having equal number of ones and minus ones, are orthogonal to the
vector~$\vec{\mathbf 1}$. They are linearly independent, and 
therefore they span $U_4$.

\begin{figure}[htbp]
\centerline{\epsffile{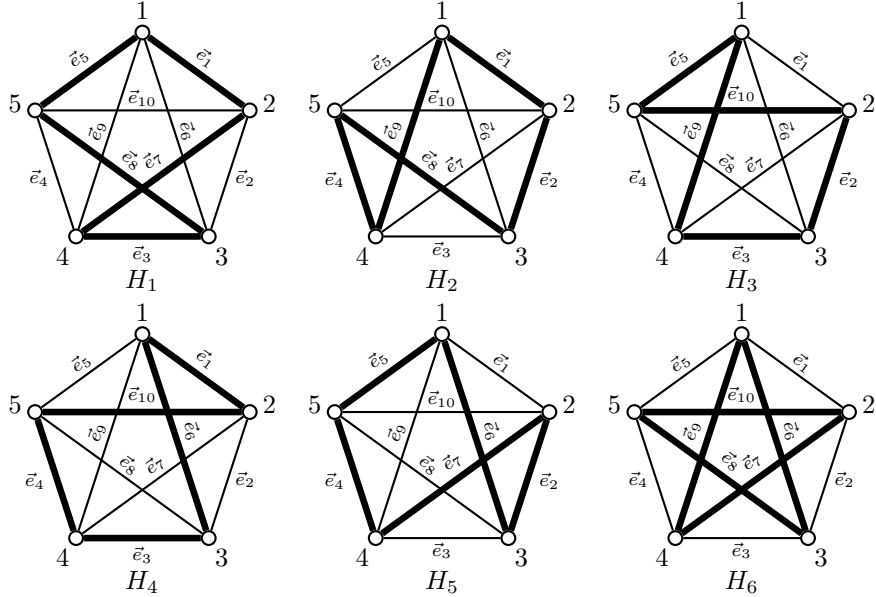}}
\caption{\small Pairs of Hamiltonian cycles in $K_5$: the ``positive''
cycles $H_i$ are drawn in bold lines, while their complements, the 
``negative'' cycles $\overline{H}_i$, are shown in thin lines.}
\label{corr3}
\end{figure}

Another collection of subsets of $E$ which is stable under the action
of\/ ${\rm S}_5$ is the set of Hamiltonian cycles $H\subset E$, that is, 
cycles that visit each vertex exactly once. A Hamiltonian cycle in $K_5$ 
can be described by a 5-cycle $c\in{\rm S}_5$ which indicates in which 
order the vertices are visited; note that $c^{-1}$ describes the same 
Hamiltonian cycle since our graph is undirected. The complement 
$\overline{H}=E\setminus H$ is also a Hamiltonian cycle which corresponds
to the permutation $c^2$ (or to its inverse $c^{-2}$). There are 24 cyclic
permutations in ${\rm S}_5$; they give rise to 12 Hamiltonian cycles 
in~$K_5$ which form 6~pairs of mutually complementary cycles: see 
Fig.~\ref{corr3}.

The vectors $\vw_k = \sum_{u\in H_k}e_u - \sum_{u\in\overline{H}_k}e_u$ 
or, more concretely,
\begin{eqnarray*}
 \vw_1 & = & (1,-1,1,-1,1,-1,1,1,-1,-1)\,, \\
 \vw_2 & = & (1,1,-1,1,-1,-1,-1,1,1,-1)\,, \\
 \vw_3 & = & (-1,1,1,-1,1,-1,-1,-1,1,1)\,, \\
 \vw_4 & = & (1,-1,1,1,-1,1,-1,-1,-1,1)\,, \\
 \vw_5 & = & (-1,1,-1,1,1,1,1,-1,-1,-1)\,, \\
 \vw_6 & = & (-1,-1,-1,-1,-1,1,1,1,1,1)\,,
\end{eqnarray*}
(once again the first five and the last five components move cyclically)
span an invariant subspace. Every edge of $K_5$ belongs to 3 ``positive''
Hamiltonian cycles and to 3 ``negative'' ones; therefore,
$\sum_{i=1}^6\vw_i=0$. It is easy to verify that the space $U_5$
spanned by these 6 vectors is in fact 5-dimensional. For every fan 
$F_i$ and for every pair $(H_j,\overline{H}_j)$, exactly two edges 
of $F_i$ belong to $H_j$, while the other two belong to $\overline{H}_j$.
Therefore, $\vv_i\perp\vw_j$ for all $i,j$, so $U_5\perp F$ where,
as before, $F=U_1\oplus U_4$.

Thus, we get a decomposition of $\Q^{10}$ into three invariant subspaces: 
\linebreak[4]
$\Q^{10}=U_1\oplus U_4\oplus U_5$. We did not prove that the subspaces
$U_4$ and $U_5$ are irreducible. The proof goes by some routine
verification using the character table of\/ ${\rm S}_5$. We omit the 
details since for our goal this fact is irrelevant: the only thing we 
wanted to show was the reducibility of the orthogonal complement 
$U_1^{\perp}=U_9$, and this statement is proved since we have shown 
that $U_9=U_4\oplus U_5$.

\bigskip

We finish this section by specifying how certain elements of ${\rm S}_5$
act on the labels of the 10 edges. By construction, the permutation
$f=(1,2,3,4,5)\in{\rm S}_5$ acts as 
$$
\varphi = (1,2,3,4,5)(6,7,8,9,10)\,.
$$
Taking a simple transposition, for example, $a=(2,5)\in{\rm S}_5$,
we get
$$
\alpha = (1,5)(2,8)(4,7)\,.
$$
Indeed, all the edges having both ends different from 2 and 5, remain
fixed, as well as the edge $\{2,5\}$ itself, while the 6 edges having
exactly one end equal to 2 or to 5 split into 3 pairs. Finally, taking
$s=(1,2)(3,5,4)\in{\rm S}_5$ we obtain
$$
\sigma = (2,5,7,6,10,9)(3,8,4)\,.
$$
Note that $s^3=(1,2)$; conjugating this element by $f$ we get all the
transpositions $(i,i+1)$ of adjacent elements. Therefore, the elements 
$s^3$ and $f$, and hence also $s$ and $f$, generate the whole group 
${\rm S}_5$. Since 
$$
\sigma\alpha\phi=1
$$ 
and the homomorphism ${\rm S}_5 \to {\rm S}_{10}$ is injective, this 
implies that the group $\langle\sigma,\alpha,\varphi\rangle$ is 
generated by $\alpha$ and $\sigma$ and is isomorphic to ${\rm S}_5$.
The action of $\langle\sigma,\alpha,\varphi\rangle\cong{\rm S}_5$ 
on the 10 edges is primitive; indeed, we could only have 2 blocks 
of 5 elements each, or 5 blocks of 2 elements each, but the presence 
of a cycle of order~6 is incompatible with the first possibility
while the presence of a single fixed point is incompatible with the
second one. The action is obviously transitive.

\section{Realization of the degree-10 action of ${\rm S}_5$ 
as the monodromy group of a Laurent polynomial}

During all this section, we systematically use various methods and
results of the theory of ``dessins d'enfants''. We will try to be
concise but clear. For all missing details the reader may address
the book~\cite{lz} (Chapters~1 and 2).

\subsection{Belyi functions and ``dessins d'enfants''}

Rational functions from $\C\P^1$ to $\C\P^1$ (and, more generally,
meromorphic functions from a Riemann surface $X$ to $\C\P^1$),
unramified outside 0, 1, and $\infty$, are called {\em Belyi functions}.
They have many remarkable properties. In particular, any such function 
$F(x)$ may be ``encoded'' in the form of a bicolored map~$M_F$ drawn
on the sphere (resp., on the surface $X$). 
Namely, let us color the points 0 and 1 in black and white respectively, 
draw the segment $[0,1]$, and define $M_F$ as the preimage 
$M_F=F^{-1}([0,1])$ of the segment $[0,1]$ with respect to the function 
$F(x)\,:\,\C\P^1\rightarrow \C\P^1$. By definition, black (resp., white) 
vertices of $M_F$ are preimages of the point $0$ (resp., of the point 1) 
and edges of $M_F$ are preimages of the segment $[0,1]$.

The segment $[0,1]$ may itself be considered as a bicolored map having
two vertices of degree~1 and a face of degree~2 containing infinity.
Clearly, $M_F$ has $n=\deg F$ edges, and the degree of a vertex $x$ 
of $M_F$ coincides with the multiplicity of $x$ with respect to $F$. 
Furthermore, each face of $M_F$ contains a pole of $F$, and twice
the multiplicity of this pole coincides with the degree of 
the corresponding face.
The map $M_F$ permits to reconstruct the monodromy group $G_F$ of $F$. 
Indeed, let $g_{0}$, $g_1$ be generators of $G_F$ corresponding to the 
loops around $0$ and $1$. Taking a base point of the covering somewhere
inside the segment $[0,1]$ we may consider that the permutations $g_0$
and $g_1$ act not on the preimages of the base point but on the preimages
of $[0,1]$, that is, on the edges of $M_F$. The permutation $g_0$
(resp., $g_1$) sends an edge $e$ to the next one in the counterclockwise
direction around the black (resp., white) vertex adjacent to $e$.
Notice that if $g_{\infty}$ is the element of $G_F$ corresponding to 
the loop around $\infty$, then $g_0g_1g_{\infty}=1$.

\begin{figure}[htbp]
\begin{center}
\epsfig{file=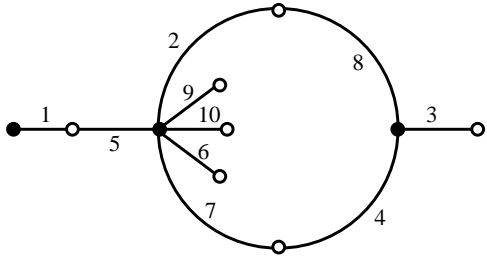,width=6.4cm}
\caption{\small Realization of ${\rm S}_5$ acting on 10 edges of a 
bicolored plane map.}\label{fig:ex}
\end{center}
\end{figure}

For example, assuming that a Belyi function $F$ corresponds to the map 
shown in Fig.~\ref{fig:ex} we may conclude that $F$ is of degree 10
(since there are 10 edges), has two poles, both of order 5 (since there
are two faces, both of degree 10), and that the corresponding 
permutations $g_0,g_1,g_{\infty}$ coincide with the permutations 
$\sigma$, $\alpha$, and $\varphi$ defined at the end of the previous 
section.

Riemann's existence theorem implies 
that  {\it for any bicolored plane map there exists a Belyi function}\/
$F(x)$ which is unique up to a composition with $x \mapsto \mu(x)$ 
where $\mu(x)$ is a  linear fractional transformation. 
In particular, since for the map shown in Fig.~\ref{fig:ex} the
permutations $g_0$, $g_1$, $g_{\infty}$ coincide with $\sigma$, $\alpha$, 
$\varphi$, this pictures ``proves'' that there exists a rational function 
$F(x)$ whose monodromy group is permutation equivalent to the action 
of ${\rm S}_5$ on 10 points discussed above.
Our next goal is to find this function explicitly.

\subsection{A system of equations for the coefficients of Belyi function, 
and its solutions}
\label{sec:belyi}

In the rest of this section we will compute a Belyi
function which produces a map isomorphic to that of Fig.~\ref{fig:ex}
as a preimage of the segment $[0,1]$. A reader not interested in the 
details of the computation may just take our word for it that the 
resulting function is the one given in (\ref{expr-for-L}), and pass 
directly to Sec.~4. We will provide not all the details but only a 
minimum allowing the reader to reproduce our results.

The black vertices of the map are the preimages of 0, or, in other words, 
they are roots of the rational function $F$ we are looking for. 
Furthermore, the vertex of degree~6 is a root of multiplicity~6, the 
vertex of degree~3 is a triple root, and the vertex of degree~1 is a 
simple root.
The freedom of choosing a linear fractional transformation $\mu(x)$ 
allows us to put these three points to any three chosen 
positions. Let us put, for example, the vertex of degree~6 to $x=0$, 
the vertex of degree~3, to $x=1$, and the vertex of degree~1, to $x=-1$. 
Then, the numerator of $F$ will take the form $x^6(x-1)^3(x+1)$.

The permutation $\varphi$ corresponds to the monodromy above $\infty$,
and it has two cycles of length~5. Therefore, the function in question 
must have two poles of degree~5, one pole inside each face of the map.
Suppose these poles to be the roots of a quadratic polynomial $x^2+ax+b$.
Then, the Belyi function in question takes the form
$$
F(x) \, = \, K \cdot \frac{x^6(x-1)^3(x+1)}{(x^2+ax+b)^5}
$$
where $K,a,b$ are constants that remain to be determined.

Here the reader may be surprised. We are looking not for an
arbitrary Belyi function but for a Laurent polynomial, aren't we? 
Then, would it not be a better idea to use the same liberty of choice 
of three parameters and to put one of the poles to $x=0$, and the other
one, to $x=\infty$? The answer is {\em no}\/: such a choice would not be a 
good idea -- at least at this stage of the computation. The reason 
is related to Galois theory and will be explained later, in 
Sec.~\ref{sec:bel-laurent}.

The white vertices of our map are the preimages of 1, or, in other 
words, the roots of the function $F(x)-1$. There are three white 
vertices of degree~2; they correspond to double roots of $F(x)-1$. 
Computing the derivative of $F$ we get
$$
F'(x) \, = \, K \cdot \frac{x^5(x-1)^2\,p(x)}{(x^2+ax+b)^6}
$$
where
$$
p(x) \, = \, (5a+2)\,x^3+(2a+10b+4)\,x^2-(a-2b)\,x-6b\,.
$$
It becomes clear that $p(x)$ is the cubic polynomial
whose roots are the three white vertices of degree~2, so the
numerator of $F(x)-1$ must have $p(x)^2$ as a factor. Note also
that the leading coefficient of this numerator is $K-1$. Thus,
we can now write down the hypothetical form of $F(x)-1$ which we
temporarily denote by $H(x)$:
$$
H(x) \, = \, \frac{K-1}{(5a+2)^2} \cdot \frac{p(x)^2\,q(x)}{(x^2+ax+b)^5}
$$
where $q(x)$ is yet unknown polynomial of degree~4, with the leading
coefficient~1, whose roots are the four white vertices of degree~1.
Denote
$$
q(x) \, = \, x^4+cx^3+dx^2+ex+f
$$
and compute the derivative of $H(x)$.

The results of the subsequent computations become more and more cumbersome. 
Their main steps go as follows. First of all, $H(x)$ is nothing else but
another representation of $F(x)-1$, so we must get in the end
$F'(x)=H'(x)$. Therefore, after having computed $H'(x)$ we ask
Maple to factor the difference $F'(x)-H'(x)$, and we get an expression
$$
F'(x)-H'(x) \, = \, {\rm Const} \cdot \frac{p(x)\,r(x)}{(x^2+ax+b)^6}
$$
where $r(x)$ is a (very huge) polynomial of degree~7. The final action 
to do is to equate $r(x)$ to zero: this means that we extract its 
coefficients and equate all of them to zero. This gives us a system
of algebraic equations on the unknown parameters $K,a,b,c,d,e,f$.

The solution of the system thus obtained using the Maple-7 package 
takes 14~seconds. It takes significantly more time to enter all the 
involved formulas and operations. And it takes even more time to find 
our way among the solutions since they are many and varied.

\subsection{Finding our way among the solutions}

\subsubsection{Maps with the same set of vertex and face degrees}

\begin{figure}[htbp]
\begin{center}
\epsfig{file=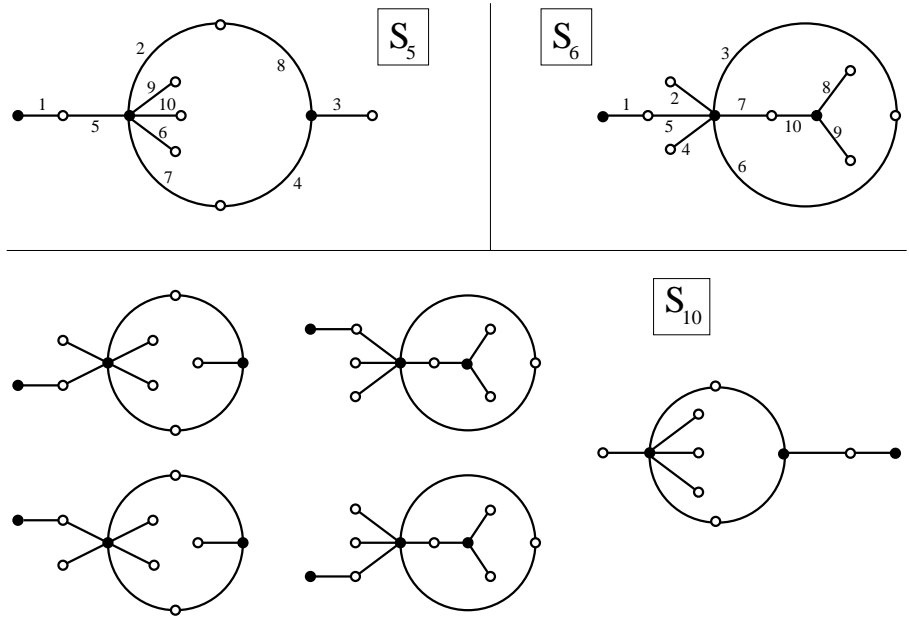,width=12cm}
\caption{\small All the seven bicolored maps with the degree partition
of the black vertices being $6^13^11^1$, that of the white vertices being
$2^31^4$, and that of the faces being $5^2$. Monodromy groups are also
indicated.}
\label{fig:vse-sem}
\end{center}
\end{figure}

If we analyse carefully the above procedure of constructing a system
of equations, we will see that the only information we have used about 
the map of Fig.~\ref{fig:ex} is the set of degrees of the black vertices, 
the white vertices, and the faces of this map. However, there exist not 
one but 7~maps having the degree partition of the black vertices equal to
$(6,3,1) = 6^13^11^1 \vdash 10$, that of white vertices equal to
$(2,2,2,1,1,1,1) = 2^31^4 \vdash 10$, and that of the faces equal to
$(5,5) = 5^2 \vdash 10$. These maps are shown in Fig.~\ref{fig:vse-sem}.
Therefore, the above computation must produce Belyi functions for 
all of them.

The picture convinces us that these 7 maps do exist. In order
to prove that there are no others we may compute the number of
triples of permutations $(g_0,g_1,g_{\infty})$ of degree~10 having 
the same cycle structure as $(\sigma,\alpha,\varphi)$ and satisfying
the equality $g_0g_1g_{\infty}=1$. For this end,
we may use, for example, the following formula due to Frobenius: 

\begin{proposition}
Let $C_1,C_2,\ldots,C_k$ be conjugacy classes in a finite group $G$.
Then the number ${\cal N}(G;C_1,C_2,\ldots,C_k)$ of $k$-tuples
$(x_1,x_2,\ldots,x_k)$ of elements $x_i \in G$ such that each
$x_i \in C_i$ and $x_1x_2\ldots x_k = 1$, is equal to
$$
{\cal N}(G;C_1,C_2,\ldots,C_k) \; = \;
\frac{|C_1|\cdot|C_2|\cdot\ldots\cdot|C_k|}{|G|} \cdot
\sum_{\chi} \frac{\chi(C_1)\chi(C_2)\ldots\chi(C_k)}{(\dim\chi)^{k-2}}\,,
$$
where the sum is taken over the set of all irreducible characters
of the group $G$. 
\end{proposition}

Applying this formula to the group $G = {\rm S}_{10}$, $k=3$, and the 
conjugacy classes $C_1$, $C_2$, $C_3$ determined by the cycle 
structures $6^13^11^1$, $2^21^4$, and $5^2$, respectively, and computing 
the irreducible characters of ${\rm S}_{10}$ using the Maple package 
{\tt combinat}, we get
$$
{\cal N}(G;C_1,C_2,C_3) \; = \; 25\,401\,600 \; = \; 7 \cdot 10!\,.
$$

None of the maps shown in Fig.~\ref{fig:vse-sem} has a non-trivial
orientation preserving automorphism; therefore, each of them admits
$10!$ different labelings.

It is useful to determine monodromy groups of the functions 
corresponding to the above maps. For the 
map in the upper left corner we know already that, by construction,
it is isomorphic to ${\rm S}_5$. For the 5 maps shown in the lower part 
of the figure, the order of the group (which can be calculated by the 
Maple package  {\tt group}, function {\tt grouporder}) is equal to 
$10!$, and therefore the group is ${\rm S}_{10}$ itself. Finally,
for the map in the upper right corner, using the same Maple package, 
or GAP, or the catalogue~\cite{ButMcK-83}, we may establish that it is 
isomorphic to~${\rm S}_6$.

\subsubsection{Galois action on maps and finding $F(x)$}\label{sec:galois}

We find the coefficients of the Belyi functions by solving a system
of algebraic equations. Therefore, there is no wonder that these
coefficients are algebraic numbers. The group of automorphisms
of the field $\overline{\Q}$ of algebraic numbers is called the 
{\em absolute Galois group}\/ and is denoted by
$\Gamma = {\rm Gal}(\overline{\Q}|\Q)$. An element of the 
group~$\Gamma$, acting simultaneously on all the coefficients of a 
given Belyi function, transforms it into another Belyi function which 
may correspond to another map.

\begin{figure}[htbp]
\begin{center}
\epsfig{file=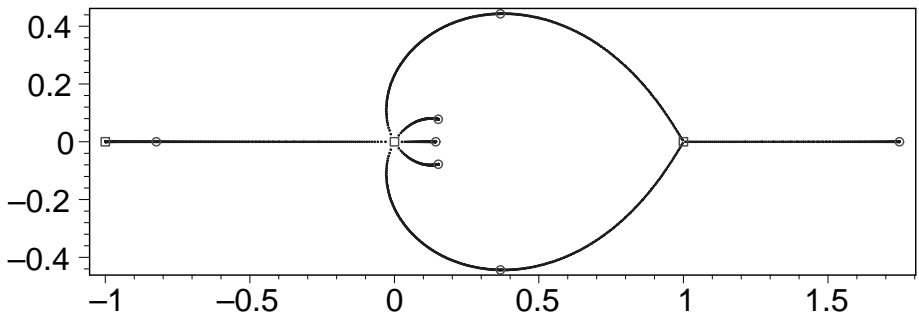,width=12cm}
\caption{\small A Maple plot of the ``dessin d'enfant'' corresponding 
to the Belyi function~$F_1$. Black vertices are marked by little squares.}
\label{fig:dess-1}
\end{center}
\end{figure}

Thus, bicolored maps split into the orbits of the above Galois action. 
The set of degrees of black and white vertices and faces is
an invariant of this action; therefore, all the orbits are finite.
Another invariant is the monodromy group. Looking once again at
Fig.~\ref{fig:vse-sem} we see that the set of 7 maps represented there
splits into at least three Galois orbits: two orbits contain each
a single element, while the set of the remaining 5 elements may constitute 
one orbit or further split into two or more orbits. The general theory 
suggests that for the singletons the coefficients of the corresponding 
Belyi functions must be rational numbers. And indeed, among our solutions
we find two such functions:
$$
F_1(x) \, = \, \frac{50000}{27} \cdot
               \frac{x^6(x-1)^3(x+1)}{(x^2+4x-1)^5}
$$
and
$$
F_2(x) \, = \, 337500 \cdot \frac{x^6(x-1)^3(x+1)}{(11x^2+4x-16)^5}\,.
$$
At this stage we simply ask Maple to draw the $F$-preimages of the 
segment $[0,1]$ and find out that the function we are looking for is~$F_1$:
just compare Fig.~\ref{fig:dess-1} with Fig.~\ref{fig:ex}. It is pictures
like that in Fig.~\ref{fig:dess-1}, obtained as Belyi preimages
of the segment $[0,1]$, which are usually called {\em dessins d'enfants}.

\bigskip

The five remaining maps constitute an orbit of degree~5 defined 
over the splitting field of the polynomial
$$
Q(t) \, = \, 85237\,t^5 - 95206\,t^4 + 
48850\,t^3 - 7456\,t^2 + 1606\,t - 226\,.
$$
This means that the coefficients of a Belyi function $F(x)$ are
expressed in terms of (more exactly, as polynomials of degree $\le 4$
in) a root of this polynomial. Taking one by one five roots we
obtain five different Belyi functions which correspond to the five
maps with the monodromy group ${\rm S}_{10}$ shown in the lower part
of Fig.~\ref{fig:vse-sem}.

Notice that besides the solutions mentioned above, our system of 
algebraic equations produces a bunch of the so-called ``parasitic 
solutions'' representing various kinds of degeneracies. Some of them are
easy to eliminate, others are not. For example, in one of
the solutions we get  $a=0$, $b=0$, which means that the denominator 
of $F$ is $x^{10}$, while its numerator contains $x^6$. This solution 
does correspond to a Belyi function, but of degree~4 instead of 10.
Another easy case is $K=1$, $a=-2/5$, which leads to a division of
zero by zero in the constant factor $(K-1)/(5a+2)^2$ of the 
function $H$ in Sec.~\ref{sec:belyi}.
More difficult cases of degeneracies also exist but we will not 
go here into further details, as well as into many 
other subtleties proper to any experimental work. The questions 
already discussed show quite well why the computation of Belyi 
functions remains a handicraft instead of being an industry.

\subsection{From a rational function to a Laurent polynomial}
\label{sec:bel-laurent}

Now we may return to the question asked in Sec.~\ref{sec:belyi} 
and explain why we decided to compute a ``generic'' Belyi function 
instead of looking from the very beginning for a Laurent polynomial.

We see that, while $F_1$ is defined over $\Q$, its two poles are not:
they are roots of the quadratic polynomial $x^2+4x-1$; concretely,
they are equal to $-2 \pm \sqrt{5}$. Any linear fractional
transformation of the variable $x$ sending one of theses poles
to~0 and the other one to $\infty$ would inevitably add $\sqrt{5}$
to the field to which belong the coefficients of Belyi functions.
Thus, the functions defined over $\Q$ would become defined over
$\Q(\sqrt{5})$, the orbit of degree~5 would become one of degree~10
(with each of the five maps being represented twice), parasitic solutions 
would also become more cumbersome (and their parasitic nature
would be more difficult to detect), and so on. And without doubt Maple 
would have a much harder work to solve the corresponding more 
complicated system of algebraic equations.

But from now on, after making all the above computations with the
simplest possible fields, we can easily transform $F_1$ into a
Laurent polynomial. The transformation
$$
z \, = \, \frac{(2-\sqrt{5})x-1}{(2+\sqrt{5})x-1}
$$
sends the pole $-2-\sqrt{5}$ to 0 and $-2+\sqrt{5}$ to $\infty$,
and also sends 0 to 1. \linebreak[4] Substituting into $F_1(x)$ 
its inverse
$$
x \, = \, \frac{z-1}{(2+\sqrt{5})z-(2-\sqrt{5})}
$$
we obtain the Laurent polynomial of Theorem~\ref{main}:
$$
L(z) \, = \, K \cdot \frac{(z-1)^6(z-a)^3(z-b)}{z^5}
$$
where
$$
K = \frac{11+5\sqrt{5}}{216}, \qquad
a = \frac{-3+\sqrt{5}}{2}, \qquad
b = \frac{7-3\sqrt{5}}{2}\,.
$$

\section{Proof of the main theorem}
\label{sec:proof}

We are looking for  Laurent polynomials $Q_j$, $0 \leq j \leq 4$, 
of the form
\be\label{form}
Q_j(z) \,=\, \sum_{k=-j}^j s_k z^{k}
\ee 
(we set $Q_0=1$) satisfying the equation (\ref{1l}). However, it is 
clear that we may multiply $Q_j$ by a constant, and also add to $Q_j$ 
an arbitrary linear combination of $Q_i$ for $i<j$, and this gives us 
another solution having the same form. Therefore, in order to achieve 
uniqueness, we impose on $Q_j$ the following three conditions:
\begin{enumerate}
\item The coefficient $s_{-j}$ is equal to one.
\item For $i=-j+1,\ldots,0$ the coefficients $s_i$ are equal to zero.
\item $\displaystyle \int_{S^1} {L}^i dQ_j = 0$ \,\, for all \,\,
      $1 \le i \le j$\,.
\end{enumerate}

The Laurent polynomial $Q_j$ has $2j+1$ coefficients; the first two 
conditions fix $j+1$ of them, while the third condition provides us 
with $j$ additional linear equations on coefficients. In order 
to ensure that the integrals in the third condition vanish, according 
to the Cauchy theorem, we must calculate the coefficients preceding
$z^{-1}$ in $L^i\cdot Q^{\prime}_j$, $1 \le j \le 4$, $1 \le i \le j$,
and set them to zero. The existence and uniqueness of solutions
will be explained later (in Step~3 of the proof). The results of the 
calculation are collected below. 

$$
Q_0 \,=\, 1\,,
$$
\vskip 0.1cm
$$
Q_1 \,=\, {\frac {{z}^{2}+1}{z}}\,,
$$ 
\vskip 0.2cm
$$
Q_2 \,=\, {\frac {-(9+4\,\sqrt{5})\,z^4+(20+8\,\sqrt{5})\,z^3+1}{{z}^{2}}}\,,
$$ 
\vskip 0.2cm
$$
Q_3 \,=\, \frac{\left(\frac {47}{2}+\frac{21}{2}\,\sqrt{5}\right) z^6 - 
\left(\frac {195}{2} + \frac {87}{2}\,\sqrt{5}\right) z^5 + 
\left(\frac{255}{2} + \frac {111}{2}\,\sqrt{5}\right) z^4 \,+\, 1}{z^3}\,,
$$ 
\vskip 0.2cm
{\small
$$
Q_4 = \frac{-(9+4\,\sqrt{5})\,z^8 + 
(130+58\,\sqrt{5})\,z^7 - 
(630+282\,\sqrt {5})\,z^6 +  
(910+406\,\sqrt {5})\,z^5 + 1}{z^4}.
$$
}

\medskip

Now, everything is ready in order to prove the main theorem. The proof
is divided into several steps.

\paragraph{Step 1.} \hspace{-2mm} First of all, we must check that the 
Laurent polynomials $Q_j$, $1\leq j \leq 4$, satisfy the equalities 
\be\la{ttyy} 
\int_{S^1} {L}^i\d Q \,=\, 0
\ee 
for {\it all}\/ $i\geq 0$ (for $Q_0$ it is obvious). For this purpose 
we may use Theorem 7.1 of \cite{pak2} and verify this equality only 
for a finite number of $i$, namely, for
\be \la{ner}
1 \,\leq\, i \,\leq\, (N - 1)\cdot\deg Q + 1,
\ee
where $N$ is the size of the orbit of the vector
$$
(1,1,1,1,1,-1,-1,-1,-1,-1)
$$
under the action of the monodromy group of $L$. In our case,
$\deg Q_j = 2j $, $1\leq j \leq 4,$ and $N = 12$; therefore, 
the maximal value of the right hand side of~\eqref{ner} is equal to 89.
The verification for all the four polynomials $Q_j$ takes less than 
one minute of work of Maple-11.

\paragraph{Step 2.} Observe that if a Laurent polynomial $Q$ is a 
solution of \eqref{ttyy} then for any polynomial $R$ the Laurent 
polynomial $\widehat Q=R(L)\cdot Q$ is also a solution of~\eqref{ttyy}.
Indeed, it is enough to prove it for $R=z^k$, $k\geq 1$. We have: 
$$
\int_{S_1} L^i \d (L^k Q) \,=\, \int_{S_1} L^{i+k} \d Q+ 
\int_{S_1} L^iQ\,\d L^k \,.
$$ 
The first integral in the right-hand side of this equality vanishes 
by \eqref{ttyy}. On the other hand, for the second integral we have: 
$$
\int_{S^1} L^iQ\d L^k = k\int_{S^1} L^{i+k-1}Q\d L =
\frac{k}{i+k}\int_{S^1} Q\d L^{k+i} = -\frac{k}{i+k}\int_{S^1} L^{k+i}\d Q\,,
$$ 
and therefore this integral also vanishes. 

\paragraph{Step 3.} The final ingredient we need is Theorem 6.7 
of~\cite{pak2} which states that if the leading degree of a Laurent 
polynomial $L$ is a prime number $p$ (in our case $p=5$), and if $Q$ 
is a {\it polynomial}\/ (that is, a common one, not a Laurent
polynomial) such that \eqref{ttyy} holds, then either 
$L(z)=L_1(z^p)$ for some Laurent polynomial $L_1$ while $Q$ is a 
linear combination of the monomials $z^{i}$ with $i$ not being 
multiples of $p$, or $Q$ is a constant. Since the Laurent polynomial 
$L$ we are working with is not of the form $L(z)=L_1(z^p)$ this result 
implies that a polynomial $Q$ cannot satisfy \eqref{ttyy} unless $Q$ 
is a constant.

This fact also explains the uniqueness of $Q_j$. Indeed, if $Q_j^{(1)}$
and $Q_j^{(2)}$ are two solutions of the equations imposed on $Q_j$
at the beginning of this section, then their difference 
$Q_j^{(1)}-Q_j^{(2)}$ is also a solution of the Laurent polynomial
moment problem. But this difference is a polynomial (since the terms
$z^{-j}$ in $Q_j^{(1)}$ and~$Q_j^{(2)}$ cancel) and therefore must
reduce to its constant term; but the constant term of this polynomial 
is equal to zero.

The uniqueness of the solution implies the non-degeneracy of the matrix
of the system, and the non-degeneracy, in its turn, implies existence.

\paragraph{Step 4.} Now, let us suppose that $Q$ is a Laurent polynomial 
satisfying \eqref{ttyy} and $m(Q)\leq 0$ is the minimal degree of a 
monomial in $Q$. Let $m(Q) =-5k_0-j_0$, where $0\leq j_0 \leq 4$ 
and $k_0\geq 0$. Then for any $c_0\in \C$ the Laurent polynomial 
$$
Q^{(1)} \,=\, Q-c_0\,L^{k_0}\cdot Q_{j_0}\,,
$$ 
is a solution \eqref{ttyy}. Furthermore, choosing an appropriate $c_0$ 
we can assume that $m(Q^{(1)})>m(Q)$ (here we use the fact that the
coefficient $s_{-j}$ in $Q_j$ is not zero). Now, if $m(Q^{(1)})=-5k_1-j_1$, 
where $0\leq j_1 \leq 4$ and $k_1\geq 0$, then, setting  
$$
Q^{(2)} \,=\, Q^{(1)}-c_1\,L^{k_1}\cdot Q_{j_1}
$$
for an appropriate $c_1$ we obtain a solution \eqref{ttyy} with 
$m(Q^{(2)})>m(Q^{(1)})$. Continuing in this way we will eventually 
arrive to a solution $Q^{(r)}$ of \eqref{ttyy} for which $m(Q^{(r)})\geq 0$. 
In view of the result cited in Step~3 such a solution should be a 
constant $c\in \C$. Therefore,
$$
Q \,=\, c \,+\, \sum_{i=0}^{r-1}c_i\,L^{k_i}\cdot Q_{j_i} \,=\,
\sum_{j=0}^4 (R_j\circ L)\cdot Q_j
$$
for some polynomials $R_0, R_1, R_2, R_3, R_4$. The theorem is proved.

\paragraph{Final remarks.} 
In general, it is not known if the reducibility of the action of the 
monodromy group $G_L$ of a Laurent polynomial $L$ of degree $n$ on the 
space $\Q^n$ always implies a non-trivial structure of solutions of the
corresponding moment problem. The only facts which follow from the 
general theory are as follows: 
\begin{itemize}
\item   The reducibility of the above action implies the existence 
        of a {\em rational function}\/~$Q$, which is not a rational function 
        in $L$, such that the generating function for the sequence 
        of the moments  
        $$
        m_i=\int_{S^1} L^i\,dQ, \ \ \ i\geq 0,  
        $$
is rational (see Sec.~8.3 of \cite{pak2}).
\item   If the above function $Q$ turns out to be a Laurent polynomial, 
        then the rationality of the generating function implies its 
        vanishing (see Theorem~3.4 of \cite{pak2}).
\end{itemize}
It would be interesting to understand in a more profound way what is 
the underlying mechanism which relates the structure of solutions of the 
moment problem for $L$ with the structure of the representation of $G_L$. 

\end{document}